\documentclass{article}

\usepackage{arxiv}

\usepackage[table]{xcolor}
\usepackage{tabularx}
\usepackage{enumitem}

\usepackage[utf8]{inputenc} 
\usepackage[T1]{fontenc}    
\usepackage{hyperref}       
\usepackage{url}            
\usepackage{booktabs}       
\usepackage{amsfonts}       
\usepackage{amsmath}
\usepackage{mathrsfs}
\usepackage{nicefrac}       
\usepackage{microtype}      
\usepackage{lipsum}
\usepackage{graphicx}
\usepackage{floatrow}
\usepackage{subcaption }
\usepackage{comment}
\usepackage[noend]{algpseudocode}
\algnewcommand{\LineComment}[1]{\State \(\triangleright\) #1}
\usepackage{algorithmicx}
\usepackage{algorithm}
\usepackage{float}
\floatstyle{plaintop}
\restylefloat{table}
\graphicspath{ {./} }

\setlist[enumerate]{
  nosep
}

\title{Learning the Markov Decision Process \\in the Sparse Gaussian Elimination}

\author{
  Yingshi Chen \thanks{Thanks to Giga for its help in this research} \\
  Giga Design Automation Co., Ltd\\
  Shenzhen 518055, China\\
  \texttt{yschen@giga.com.cn} 
}

\begin{document}
\maketitle 
\begin{abstract}
 We propose a learning-based approach for the sparse Gaussian Elimination. There are many hard combinatorial optimization problems in modern sparse solver. These NP-hard problems could be handled in the framework of Markov Decision Process, especially the Q-Learning technique. We proposed some Q-Learning algorithms for the main modules of sparse solver: minimum degree ordering, task scheduling and adaptive pivoting. Finally, we recast the sparse solver into the framework of Q-Learning.
 
 Our study is the first step to connect these two classical mathematical models: Gaussian Elimination and Markov Decision Process. Our learning-based algorithm could help improve the performance of sparse solver, which has been verified in some numerical experiments.
\end{abstract}

\keywords{Gaussian Elimination \and Markov Decision Process \and Sparse Solver \and Q-Learning \and Matrix Computation \and Task Scheduling}

\section{Introduction}
Gaussian elimination(GE) is one of the oldest and most important numerical method\cite{1996Matrix,grcar2011mathematicians,yuan2012jiu}. It's the most widely used method for solving linear systems:
\begin{equation}
    Ax=b
\end{equation}
where $A$ is a square matrix with rank $n$, $b$ is the use specified right hand side and $x$ is the solution. In many literatures, this method is also called direct method \cite{duff2017direct,demmel200010}, or direct sparse solver in the case of sparse matrix.

At the age of Gauss (1777–1855), the number of equations is very small. In recent decades, the rank of sparse matrix is getting larger and larger, from thousands to millions and billions and trillions. Its solution becomes the time and resource bottleneck of many applications. 

In addition to the explosion of scale, another main challenge is how making efficient use of the hierarchical architecture of modern computers\cite{dongarra1998numerical}. Only a few highly optimized dense matrix operators (such as BLAS\cite{dongarra1990set,anderson1999lapack,kaagstrom1998gemm}) could make full use of the great performance of cores(CPU,GPU,TPU...). These dense block operation is hard to implement in sparse solver. How construct these blocks? How scheduling these blocks? How reduce the communication in these blocks? These all involve complex combinatorial optimization problems. Most problems are actually NP-hard \cite{bovet1994introduction}, which cannot be solved in polynomial time. There are more problems unique to sparse solver, such as fill-in ordering (Yannakakis proved it is NP-complete\cite{yannakakis1981computing}). Therefore, sparse solver is not just a simple extension of GE, but do need a new framework to handle these problems.

Modern sparse solver \cite{schenk2004solving,li2003superlu_dist,amestoy2000mumps,Chen2018} is based on the operators on blocks. The first version is the frontal solver by Bruce Irons \cite{irons1970frontal}. It builds a LU or Cholesky decomposition from the assembly of element matrices (fronts). Duff and Reid\cite{duff2017direct} improved this idea to multifrontal solver, which deals with several independent fronts at the same time\cite{liu1992multifrontal}. Their research has laid the foundation of modern sparse solver. And now we are still walking along the road they pointed out. Maybe this paper is the first attempt to improve multifrontal method from a more general framework, a learning-based framework.

\subsection{Main Modules of Sparse Solver} 
Modern sparse solver could be divided into three stages (there are also four stages division, that is, an independent fill-in order stage):
\begin{enumerate}
\item Symbolic analysis: includes many tasks on the symbolic (non-zero) pattern of matrix.       
\begin{itemize}
    \item Transform a matrix into a specific structure (For example, block triangular form).
    \item Reorder the matrix to reduce fill-in (Detail in section \ref{sec:md}) 
    \\The minimum degree order algorithm is NP-hard\cite{yannakakis1981computing}.
    \item Get the pattern of the L and U matrix. Construct frontals from the pattern.
    \\There are some greedy algorithm to construct frontals. But the optimal frontals partition is NP-hard. We would prove it in a later paper.
    \item Static parallel task scheduling on these frontals (Detail in section \ref{sec:mdp_dag}).
    \\The static task scheduling problem is NP-hard \cite{graham1969bounds}.
    \item More, such as the elimination tree \cite{gilbert1993elimination,liu1990role}.
\end{itemize}
Finally, this stage would get a permutation matrix $P$.

\item Numeric factorization: computes the $L$ and $U$ factors in frontal format.     \\
\begin{equation}
LU=P^{T}AP
\end{equation}
 Pivoting is an important technique to improve numerical stability. For more details of adaptive pivoting, see section \ref{sec:mdp_pitvoting}.
 
This is usually the most time-consuming stage. The key to reducing time and improving efficiency is dynamic (online) task scheduling. It's NP-hard \cite{graham1969bounds} and much more complicated than static scheduling. For more details, see section \ref{sec:mdp_dag}.

\item Solve:  performs forward and backward substitution to get solution $x$.   

\begin{equation}
x=\left (P^T  \right )^{-1}\left (U  \right )^{-1} L^{-1} P^{-1} b
\end{equation}
Dynamic task scheduling is also needed to improve the effeciency.
\end{enumerate}

\subsection{Learn the Sparse Solver by Markov Decision Process}
It's clear that there are many combinatorial optimization problems in sparse solver. To handle these problems, Markov Decision Process (MDP in short) \cite{bellman1957markovian,gosavi2015solving} is a smart and powerful model.
We would use several sections to reveal this seemingly surprising discovery, to reveal many MDP in sparse solver. And how to learn this process. And how to improve the performance from learning-based method.

These methods have recently been in great glory in academia and industry. For example, Google's AlphaGO\cite{silver2016mastering} defeated the go master Lee Sedol. The engine of AlphaGO is based on deep reinforcement learning, whose theoretical basis is MDP. It's time to use these powerful modern weapons to tackle old Gaussian elimination problem.

\section{Background}\label{sec:Background}
In this section, we give some background knowledge which are needed in the following analysis and derivation.

\subsection{Gaussian Elimination with Pivoting}    \label{back:GE}
Let's first  review the classical Gaussian Elimination with Pivoting. For more detail, please see \cite{1996Matrix,duff2017direct,foster1997growth,sherman1978algorithms}.

\begin{algorithm}[H]
\caption{Gaussian Elimination with Pivoting}    \label{alg:GE}
\begin{flushleft}
\textbf{Input:} $A\in M_{n}$ 
\end{flushleft}

\begin{algorithmic}[1]
\State $A_1=A$
\For{$k=1,2,\cdots,n-1$}
    \State Apply pivoting strategy (partial pivoting\cite{sherman1978algorithms}, rook pivoting\cite{foster1997growth},...)
    \State Select pivot element ant its position is $(r,c)$
    \If { $\left ( r,c \right ) \neq \left ( k,k \right )$ } 
    \State Interchange rows $A_k[r,:] \leftrightarrow A_k[k,:]$ and columns $ A_k[:,c] \leftrightarrow A_k[:,k]$
    \EndIf
    \LineComment{Update the un-eliminated part of the matrix}
    \For{$i=k+1,\cdots,n$}
    \State  $A_k[i,k]=A_k[i,k]/A_k[k,k]$
        \For{$j=k+1,\cdots,n$}
            \State  $A_k[i,j]=A_k[i,j] - A_k[i,k] \times A_k[k,j] $
        \EndFor
    \EndFor
\EndFor

\end{algorithmic}
\end{algorithm}

As algorithm \ref{alg:GE} shows, at step $k$, the value and pattern of $A_{k}$ depends on $A_{k-1}$ and the pivot position. Different pivoting strategies will lead to different  $A_{k}$. For a input matrix $A$, we define all possible $\left \{ A_k | k=1,2,\cdots, n \right \}$ as the elimination set of $A$, denoted by $E_{n}(A)$, or simplified as $E_{n}$ in this paper. 

The following figure 1 shows the elimination process of a simple $3 \times 3$ matrix. The first step is column partial pivoting to swap rows $R_3$ and $R_1$. And in this process, $E_{n}(A)=\left \{A_0,A_1,A_2,A_3\right \}$.
\begin{figure}[!ht]
\includegraphics[width=\textwidth]{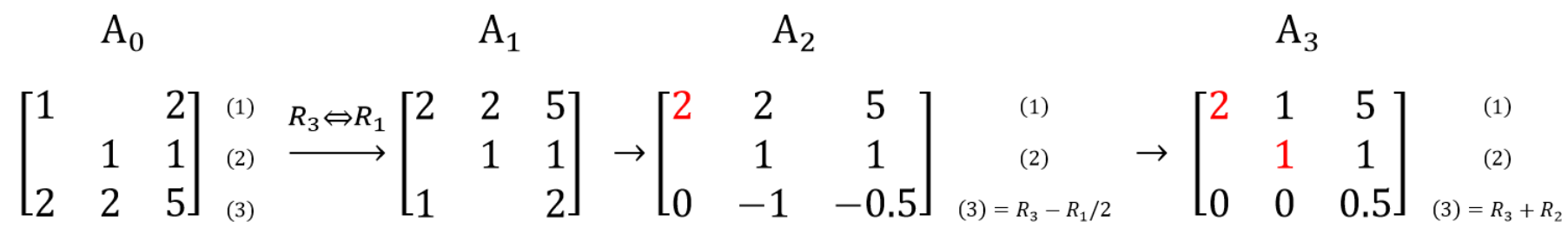}
\caption[position=bottom]{ Gaussian elimination process of a $3 \times 3$ matrix}
\end{figure}

\subsection{BLAS and Frontal}    \label{back:blas}
BLAS is an acronym for Basic Linear Algebra Subprograms\cite{dongarra1990set,anderson1999lapack,kaagstrom1998gemm}. It is the cornerstone of modern high performance computing. BLAS-1 routines perform scalar, vector and vector-vector operations, BLAS-2 routines perform matrix-vector operations, and BLAS-3 routines perform matrix-matrix operations. Many operations in frontal is BLAS-3 routines.

\subsection{Distributed/Heterogeneous Computing Systems} \label{back:distributed}
Today's distributed computing systems often contain various processors (CPU, GPU, TPU...). So it's also Heterogeneous Computing\cite{karatza2002load}. This complex system has amazing peak computing power. But the question is how to exert its power. 
There are many problems would reduce the efficiency of the whole system: 
\begin{itemize}
\item load imbalance among processes; 
\item overhead during communication; 
\item synchronizations for computations; 
\item online scheduling. A widely used model is the following DAG model.
\end{itemize}

\subsection{DAG(directed graph with no directed cycles) based Task Scheduling}    \label{sec:DAG}
We could use a task graph to  describe all tasks in sparse solver, where the nodes correspond to the tasks on frontals, and the edges represent dependencies between tasks. This DAG model has been widely used in task scheduling. \cite{wu2015hierarchical} use a hierarchical DAG to schedule tasks between CPU and GPU, whereby the CPU has small tiles and the GPU has large tiles. \cite{gupta2002improved} not only uses task-DAG, but also uses data-DAG to deal with the data-dependency in sparse GE.
The DAGuE framework in \cite{bosilca2012dague} implement a complex DAG based scheduling engine, which could deal with many low level architectural issues such as load balancing, memory distribution, cache reuse and memory locality on non-uniform memory access (NUMA) architectures, and communications/computations overlapping. For more detail of task scheduling, please see section \ref{sec:mdp_dag}

\subsection{Markov Property}    \label{sec:Markov}
Markov property is named after the brilliant Russian mathematician Andrey Markov. He said: in some process, the future does not depend on the past, only depends on the present state! 
This is a very interesting and powerful property. Whether you feel it or not, many greedy/progressive algorithms are based on this property. Especially many algorithms in sparse solver. For example, MD (minimum degree order) always selects node with minimum degree at current state. Partial Pivoting always select a element in the current column/row. In the process of GE, the value and degree of each element are varying from step and step. But we don't care how they changed in the past, we only check the value at the current step. This also exists in the algorithm of task scheduling. Markov property is a prerequisite for many models, such as MDP mentioned below.

\subsection{Markov Decision Process (MDP)}    \label{sec:MDP}
The Markov decision process (MDP in short) \cite{bellman1957markovian,gosavi2015solving} is a smart and powerful model. As its name suggests, this model could find the solution (decision process) on the Markov property. The standard MDP includes four elements $\left ( E_n,O,P_o,R_o \right )$\cite{bellman1957markovian}:
\begin{itemize}
\item $E_n$ is the state space of MDP.
\item $O$ is all actions. We also use $O_s$ to denote the set of actions available from state $s$,
\item ${P_{o}(s,s')=\Pr(s_{t+1}=s'\mid s_{t}=s,o_{t}=o)}$ is the probability that action $o$ in state $s$ at time $t$,
\item $R_{o}(s,s')$ is the immediate reward (or expected immediate reward) received after transitioning from state $s$ to state $s'$, due to action $o$.
\end{itemize}

A transition matrix $M$ is a square matrix used to describe the transitions of state space. 

\begin{equation}
M=
\begin{bmatrix} 
	M_{1,1} & \cdots & M_{1,j} & \cdots & M_{1,n} \\
	& & \vdots   \\
	M_{i,1} & \cdots & M_{i,j} & \cdots & M_{i,n} \\
	& & \vdots   \\
	M_{n,1} & \cdots & M_{n,j} & \cdots & M_{n,n} \\
\end{bmatrix}
\end{equation}

Each element $ M_{i,j} $ is a non-negative real number representing a probability $ Pr(j|i) $ of moving from state $i$ to state $j$.
Each row summing of $M$ is 1: $ \sum_{j} M_{i,j}=1 $.

MDP is the theoretical foundation of Reinforcement Learning\cite{kaelbling1996reinforcement}. A widely-used Reinforcement Learning technique is Q-Learning.

\subsection{Q-Learning}             \label{sec:Q}
The MDP model presents a general paradigm and an abstract mathematical framework. For practical problems, there are many powerful techniques, such as Q-Learning\cite{watkins1989learning,watkins1992q}. "Q" refers to the expected rewards for an action taken in a given state\cite{watkins1989learning}. At each step, Q-Learning would take some action on the estimate of Q-values. If we know "Q" value of any action in a state, then it's simple to find the optimal solution of MDP. 

The following is a typical framework of Q-Learning\cite{sutton2018reinforcement}.
\begin{algorithm}[H]
\caption{A general framework of Q-Learning}
\begin{flushleft}
\textbf{Init Parameters:} 
Learning rate $\alpha \in \left ( 0,1 \right ]$, discount rate $\gamma $ and reward function $R$
\end{flushleft}

\begin{algorithmic} [1]
\State Initialize $Q\left ( s,a \right )$ for all state $s$ and action $a$. For the terminal state, $Q(terminal,\cdot)=0$
\For{ each episode}
    \State Pick a initial state $s$
    \While {true}
        \State Choose action $a$ for the current state $s$, using policy derived from $Q$ 
        \State Take action $a$ from some policy $\pi$ then enter the next state $s'$
        \State Get reward $r = R\left ( s,a \right )$
        \State Update Q by some method, for example:
        \State \quad $Q\left ( s,a \right ) \leftarrow Q\left ( s,a \right )+\alpha \left [ r+\gamma Q\left ( s',a' \right )-Q\left ( s,a \right )  \right ]$
        \State Set $s=s'$
        \If {$s$ is terminal}
            \State break
        \EndIf
    \EndWhile
\EndFor
\end{algorithmic}
\end{algorithm}
This framework shows a key advantage of Q-Learning - avoiding the usage of transition matrix. In the case of sparse solver, the probability in the transition matrix is hard to estimate. Or the transition matrix needs huge memory
The number of states of a million-order matrix is astronomical!). Instead, we could always get the reward or Q-value. So in practice, Q-Learning is more suitable for studying various combinatorial problems in GE. Watkins\cite{watkins1992q} prove that it would converge to the optimum Q-values with probability 1.So it's a reliable technique for the combinatorial optimization problems in sparse solver. 

\subsection{Offline Q-Learning} \label{sec:offline}
Solution time is one of the most important metrics of sparse solver. We could avoid the training time in Q-Learning by offline technique. That is, first collect many matrices and train the model on these sample matrices. Then apply the learned Q-value to solve a matrix. This offline Q-Learning (or batch Q-Learning) has recently regained popularity as a viable path towards effective real-world application. For example, Conservative Q-learning\cite{kumar2020conservative,levine2020offline}.

\begin{algorithm}[H]
\caption{Conservative Q-learning}
\begin{flushleft}
\textbf{Init Parameters:} Learning rate $\alpha \in \left ( 0,1 \right ]$, discount rate $\gamma $ and rewards function $R$
\item[] \textbf{Training stage:} 
\\ \hspace*{0.2in} Learn $Q\left ( s,a \right )$ using offline data-set
\item[] \item[] \textbf{Inference stage:}
\end{flushleft}

\begin{algorithmic} [1]
    \State Pick a initial state $s$
    \While {true}
        \State Choose action $a$ for the current state $s$, using policy derived from $Q$ 
        \State Take action $a$ from learned policy $\pi$ then enter the next state $s'$
        \State Get reward $r = R\left ( s,a \right )$
        \State Set $s=s'$
        \If {$s$ is terminal}
            \State break
        \EndIf
    \EndWhile

\end{algorithmic}
\end{algorithm}

\section{Markov Decision Process in Sparse Solver} \label{sec:mdp_ss}
In this section, we list many combinatorial optimization problems in sparse solver, which all have Markov Property(\ref{sec:Markov}). And there is implicit or explicit reward for the decision of each step. So all these algorithms could be regarded as a Markov Decision Process (MDP). We would list all the state space, actions and rewards in these algorithms. Next, we further use Q-learning to improve these algorithms.

\subsection{MDP in Minimum Degree Order} \label{sec:md}
Minimum degree (MD) order\cite{duff1974comparison,george1977application,eisenstat1977yale,george1980fast,amestoy1996approximate} is an important technique to reduce the fill-in. Its principle is simple and rough, as shown in Algorithm \ref{alg:MD}. At each step, it will always pick the node with minimum degree (may be approximate or weighted value). In the practical implementation, the key is how to update degree quickly without reducing the quality of sorting. For example, the State-Of-The-Art AMD\cite{amestoy1996approximate} uses an approximate formula on the upper-bound of degree. Surprisingly, this approximation can even get less fill-in than those methods on the accurate degree in many problems. This seemingly 'unreasonable' shows the difficulty of combinatorial optimization problems again.
\begin{algorithm}[H] 
\caption{The minimum degree order algorithm}      \label{alg:MD}
\begin{flushleft}
\textbf{Input:} 
Construct a graph $G(V^{0},E^{0})$ from the nonzero pattern of $A_0$
\end{flushleft}

\begin{algorithmic}[1]
\State For each node $i\in V^{0}$, set $d_i$ to be the degree of $i$
    \State $k=1$
    \While { $k \leq n$ }
        \State Select node $p\in V^{k-1}$ that minimizes $d_p$. 
        \State Update permuting vector   $P[k] = p$
        \State Eliminate $p$, Update $G(V^{k-1},E^{k-1})$ to $G(V^{k},E^{k})$
        \State For each node $i$ adjacent to $p$, update $d_i$
        \State $k=k+1$
    \EndWhile
\Return The permuting vector $P$
    
\end{algorithmic}
\end{algorithm}

As algorithm \ref{alg:MD} shows, the degrees of nodes are varying from step and step. But we don’t care how they changed in the past, we only check the degree at the current step. That's Markov Property, decisions made does not depend on the past, only depends on the present state! 

\subsubsection{Action and Reward in Minimum Degree Order} 
At each step of MD method, we would always pick a node to reduce the fill-in. So the negative number of fill-in could be regarded as the reward value. And different strategies correspond to different actions. In addition to AMD mentioned above, there are more choices:
\begin{itemize}
\item MMDF - modified minimum deficiency\cite{ng1999performance}.
\item MMMD - modified multiple minimum degree\cite{ng1999performance}.
\item MIND - Minimum Increase In Neighborhood Degree \cite{ng1999performance}. 
\item AMIND - Approximate Minimum Increase In Neighborhood Degree\cite{rothberg1998node}.
\item MMF - Minimum Mean Local Fill \cite{ng1999performance}.
\item AMMF - Approximate Minimum Mean Local Fill\cite{rothberg1998node}.
\end{itemize}
So the action space would be $\mathscr{A}_{md} = \left \{ \textrm{AMD}, \textrm{MMDF}, \textrm{MMMD}, \textrm{MIND}, \textrm{AMIND}, \textrm{MMF}, \textrm{AMMF}... \right \}$. And the reward $R_{md} = -N_{fill-in} $.

\subsubsection{Adaptive Minimum Degree Order by Q-Learning} 
In all previous research\cite{duff1974comparison,george1977application,eisenstat1977yale,george1980fast,amestoy1996approximate}, only one method used in MD process. No research tries to use different strategy in different step. Since all the metric in $\left \{ \textrm{AMD}, \textrm{MMDF}, \textrm{MMMD}, \textrm{MIND}, \textrm{AMIND}, \textrm{MMF}, \textrm{AMMF}... \right \}$ are reasonable, why not choose different strategies according to different patterns? This choice can be learned through offline Q-Learning algorithm(section \ref {sec:offline}). So we propose a Q-Learning based MD order method, as Algorithm \ref{alg:MD-Q} shows:
\begin{algorithm}[H]
\caption{An offline Q-Learning framework for minimum degree algorithm}      \label{alg:MD-Q}
\begin{flushleft}
\textbf{Training stage:} 
    \\ \hspace*{0.2in} Learn $Q\left ( s,a \right )$ on the MD process of many matrices.  
    \\ \hspace*{0.2in} $Q$ table restores the choice of strategy at each state.
\\ \textbf{Inference stage:}. 
\end{flushleft}

\begin{algorithmic} [1]
    \State The initial state $s$ is just $A_1=A$(input matrix)
    \While { $k \leq n$ }
        \State For the current state $A_k$, take action $a$ from learned $Q$ table. Then select node $p\in V^{k-1}$
        \State Eliminate $p$
        \State Get reward $r = R\left ( A_k,p \right )$
        \State $k=k+1$
    \EndWhile
\end{algorithmic}
\end{algorithm}

\subsection{MDP in Task Scheduling}    \label{sec:mdp_dag}
Task scheduling is the key to get high performance, especially in the distributed and heterogeneous computing\cite{sinnen2007task,luo2021learning,glaubius2012real}. The main target (objective function, reward) is to minimize the solution time (or makespan) and resource consumption (memory, power, ...). A specific goal for sparse solver is the precision. We would discuss it in a later paper. And for most problem, the direct solver would always give a solution with reasonable precision as long as the solution is completed.

This problem has been deeply and widely studied.  In recent years, the research on MDP based scheduling is very active and there are many good ideas \cite{sinnen2007task,luo2021learning,glaubius2012real,shyalika2020reinforcement}. Unfortunately, these ideas are not applied to sparse solver. To our knowledge, in the case of matrix computation, only a few papers \cite{agullo2015bridging,grinsztajn2020geometric,grinsztajn2021readys} have relevant research. They successfully apply MDP based algorithm to dense Cholesky decomposition. In the work of \cite{grinsztajn2020geometric,grinsztajn2021readys}, their Q-Learning method has the following module: 1) State space is actually from the embedding of each task to graph neural network\cite{kipf2016semi}. 
2)  Reward function is final makespan given by the whole scheduling trajectory, normalized by a baseline duration. 
3)  Action space consists in selecting an available task or in doing nothing (pass).
They trained this model with A2C\cite{mnih2016asynchronous} (a policy gradient method).  

For sparse solver, there is still no similar study. There are great difference between dense and sparse solver. They have different state space, action space and reward function. There are lots of work needed to apply MDP based task scheduling. 

\subsubsection{Action and Reward in Task Scheduling} 
There are many goals, as listed below. So many goals show the importance of scheduling.
\begin{itemize}
\item TIME - the wall time spent by the solver
\item MEMORY - memory consumption \cite{lacoste2015scheduling} 
\item BALANCE - workload balance balancing between the resources
\item OVERLAP - Efficient Communication/Computation Overlap \cite{marjanovic2010overlapping,sao2018communication}
\item POWER - Reduce the power and energy consumed\cite{aguilar2020performance}.
\item LOCALITY - 
\end{itemize}
So the action space would be $\mathscr{A}_{scheduling} = \left \{ \textrm{TIME}, \textrm{MEM}, \textrm{BALANCE}, \textrm{OVERLAP}, \textrm{POWER}, \textrm{LOCALITY}... \right \}$.
And the reward is a  weighted sum of these metrics $R_{scheduling} = \alpha \textrm{TIME} + \beta \textrm{MEM} + \gamma \textrm{BALANCE} + \cdots $.

\subsubsection{Improve Task Scheduling by Q-Learning}
For the sparse solver, the task scheduling is even more important and we propose the following algorithm:
\begin{algorithm}[H]
\caption{An Offline Q-Learning framework for the task scheduling of sparse solver} \label{alg:schedule-Q}
\begin{flushleft}
\textbf{Training stage:} 
    \\ \hspace*{0.2in} Learn $Q\left ( s,a \right )$ on the task scheduling of many matrices
    \\ \hspace*{0.2in} 
\\ \textbf{Inference stage:}
\end{flushleft}

\begin{algorithmic} [1]
    \State The initial state $s$ is just $A_1=A$(input matrix)
    \While { $k \leq n$ }
        \State For the current state $A_k$, take action $a$ from learned $Q$ table. 
        \State Run one or more task, Update task DAG.
        \State Get reward $r = R\left ( A_k,p \right )$
        \State $k=k+1$
    \EndWhile
\end{algorithmic}
\end{algorithm}

\subsection{Action and Reward in Adaptive Pivoting Method}       \label{sec:mdp_pitvoting}

Pivoting is an important technique to improve numerical stability. The partial pivoting listed in algorithm 1 is only one choice. There are more pivoting techniques. At each step of GE, there are many choices. So we could use learning-based method to pick one pivot action. That's the basic principle of adaptive pivoting. Here are some common pivoting methods:
\begin{itemize}
\item PP - Partial pivoting\cite{sherman1978algorithms}
\item RP - Rook pivoting\cite{foster1997growth}
\item CP - Complete pivoting\cite{edelman1992complete}
\item SPP - Supernodal Partial-Pivoting\cite{demmel1999supernodal}
In the case of multifrontal solver, the pivoting could be applied only in the current frontal(supernodal). But in some case, we still need to find proper pivoting element in all uneliminated frontals. 
\item SKIP - No pivoting(For example, there is no need to do pivoting in the Cholesky factorization)
\\supernodal partial pivoting version
\item RBT - Recursive butterfly transforms \cite{lindquist2020replacing}
\end{itemize}
So the action space of pivoting is $\mathscr{A}_{pivoting} = \left \{ \textrm{PP}, RP, CP, SPP, SKIP, RBT... \right \}$. 

An reasonable and classic reward value is the negative growth factor $\rho_n$\cite{wilkinson1988algebraic,foster1997growth,higham1989large}. 
\begin{equation}
R_{pivoting} = -\rho_n = -\frac{max_{i,j,k}\left |a_{i,j}^k  \right |}{max_{i,j}\left |a_{i,j}  \right |}\geq 1
\end{equation}
That is, pivoting method should reduce the growth factor as much as possible. As the outstanding Wilkinson pointed 60 years ago\cite{wilkinson1988algebraic}, low value of $\rho_n$ means high numerical stability, that's the goal of pivoting.

\section{Recasting and Improve Sparse Solver by Q-Learning}
In the previous section, we analyzed most algorithms in the sparse solver and found that these algorithms can be described by MDP. We further proposed how to use Q-Learning technology to improve these algorithms. Based on these algorithm modules, we propose a novel offline Q-Learning framework for sparse Gaussian Elimination. This framework has three basic modules $\left ( E_n,\mathscr{A},Q \right )$:
\begin{enumerate}
\item State space $E_n$ - the elimination set includes all possible eliminated matrix $A_n$(section \ref{back:GE}).
\item Action space $\mathscr{A}$ - table \ref{Table:actions} lists all actions in different Markov Decision Process(MDP). It also lists the corresponding reward metrics. 
\item $Q$ is the expected rewards for an action taken in a given state $s$.
\end{enumerate}

\begin{table}[ht]
\begin{center}
\caption{Action Space and Rewards of MDP in Sparse Solver } \label{Table:actions}
\begin{tabularx}{\textwidth}[t]{XX}

\arrayrulecolor{black}\hline    \arrayrulecolor{green}\hline
\textbf{{MDP 1 - Minimum Fill-in}} & \\
\hline
Reward $ = -N_{fill-in} $ & 
\begin{minipage}[t]{\linewidth}%
Action Space = \{ AMD, MMDF, MMMD, \\ MIND, AMIND, MMF, AMMF... \}
\end{minipage}\\

\arrayrulecolor{black}\hline    \arrayrulecolor{green}\hline
\textbf{{MDP 2 - Task scheduling}} & \\
\hline
Reward = $\alpha \textrm{TIME} + \beta \textrm{MEM} + \gamma \textrm{BALANCE} + \cdots $&
\begin{minipage}[t]{\linewidth}%
Action Space = \{ TIME, MEM, BALANCE, OVERLAP, POWER, LOCALITY... \}
\end{minipage}\\

\arrayrulecolor{black}\hline    \arrayrulecolor{green}\hline
\textbf{{MDP 3 - Pivot for numerical stability}} & \\
\hline
Reward $ = -\rho_n = -\frac{max_{i,j,k}\left |a_{i,j}^k  \right |}{max_{i,j}\left |a_{i,j}  \right |}$ &
\begin{minipage}[t]{\linewidth}%
Action Space = \{ PP, RP, CP, SPP, SKIP, RBT ... \}
\end{minipage}\\

\arrayrulecolor{green}\hline        \arrayrulecolor{black}\hline    
\end{tabularx}
\end{center}
\end{table}

Based on the algorithm \ref{alg:GE} ("Gaussian Elimination with Pivoting"), algorithm \ref{alg:MD-Q} ("An offline Q-Learning framework for minimum degree algorithm") and algorithm \ref{alg:schedule-Q} ("An offline Q-Learning framework for the task scheduling of sparse solver"), we would finally unify these algorithms, and recasting the sparse GE in the following algorithm  \ref{alg:GAUSS-Q} - "Sparse Gaussian Elimination in the framework of Q-Learning".

\algnewcommand\algorithmicswitch{\textbf{switch}}
\algnewcommand\algorithmiccase{\textbf{case}}
\algnewcommand\algorithmicassert{\texttt{assert}}
\algnewcommand\Assert[1]{\State \algorithmicassert(#1)}%
\algdef{SE}[SWITCH]{Switch}{EndSwitch}[1]{\algorithmicswitch\ #1\ \algorithmicdo}{\algorithmicend\ \algorithmicswitch}%
\algdef{SE}[CASE]{Case}{EndCase}[1]{\algorithmiccase\ #1}{\algorithmicend\ \algorithmiccase}%
\algtext*{EndSwitch}%
\algtext*{EndCase}%

\begin{algorithm}[H]
\caption{Sparse Gaussian Elimination in the framework of Q-Learning} \label{alg:GAUSS-Q} 
\begin{flushleft}
\textbf{Init Parameters:} Learning rate $\alpha \in \left ( 0,1 \right ]$, discount rate $\gamma $ and rewards function $R$
\item[] \textbf{Pre-Training  with many sample matrices:} 
\\ \hspace*{0.2in} Learn $Q\left ( s,a \right )$ using offline data-set
\item[] \item[] \textbf{Solve a matrix $A\in M_{n}$:}
\end{flushleft}

\begin{algorithmic} [1]
\State $A_1=A$
\For{$k=1,2,\cdots,n-1$}
    \State For the current state $A_k$, select action $a$ (maybe order, pivoting, task scheduling...) on the Q value.
    \Switch{action $a$}
    \Case{Minimum Degree Order:}
      \State{Update matrix pattern, ...}
    \EndCase
    \Case{Pivoting:}
        \State Get pivot element ant its position is (r,c)
        \If { $\left ( r,c \right ) \neq \left ( k,k \right )$ } 
        \State Interchange rows $A_k[r,:] \leftrightarrow A_k[k,:]$ and columns $ A_k[:,c] \leftrightarrow A_k[:,k]$
        \EndIf
        \State{Update the un-eliminated part of the matrix}
    \EndCase
    \Case{Scheduling:}
      \State{Update task DAG, ...}
    \EndCase
    \EndSwitch

\EndFor
\Return
\end{algorithmic}
\end{algorithm}

As the algorithm \ref{alg:GAUSS-Q} shows, we would first learn Q-vale in the offline training stage. These learned Q-values would help to find  a good "policy" at each step. The convergence of this algorithm is proved in \cite{watkins1992q,melo2001convergence}. So it would find an optimal policy to maximizing the total reward over any and all successive steps. We listed common reward values in table 1. The process of maximizing reward corresponds to an accurate and efficient GE process: fewer fill-in, high numerical stability, good parallelism, low communication overhead, ...

Now, we have linked two seemingly unrelated areas Q-Learning and Sparse Solver. This is not a coincidence or a blunt explanation. Most problems in sparse solvers are combinatorial optimization. And MDP is a powerful framework for these NP-hard problems. In addition to the Q-learning introduced in this paper, there are more powerful reinforcement learning methods. The reinforcement learning methods are special case of MDP and can be easily integrated into the algorithm\cite{shyalika2020reinforcement,kaelbling1996reinforcement,grinsztajn2021readys}. We would expect more improvements from these tools and try more tools \cite{chen2020deep} in later study. 

\subsection{A demo case} 
\begin{figure}[H]
\begin{subfigure}{1.0\textwidth}
  \centering
  \includegraphics[width=\textwidth]{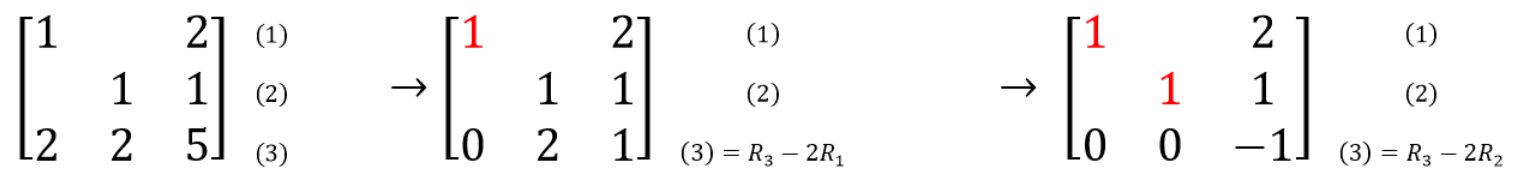}
  \caption{Minimize computation}
  \label{fig:sfig1}
\end{subfigure}%

\begin{subfigure}{1.0\textwidth}
  \centering
  \includegraphics[width=\textwidth]{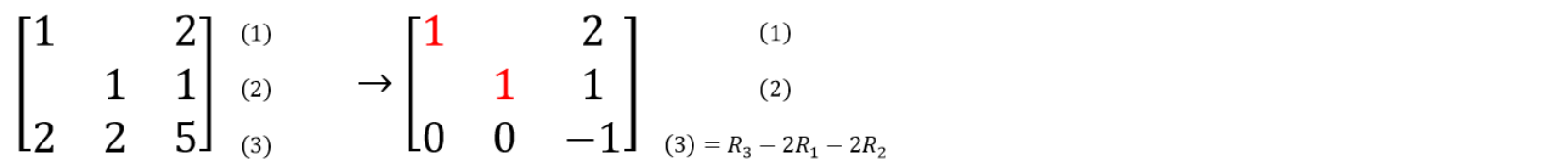}
  \caption{High parallelism}
  \label{fig:sfig2}
\end{subfigure}

\begin{subfigure}{1.0\textwidth}
  \centering
  \includegraphics[width=\textwidth]{F1}
  \caption{High numerical stability}
  \label{fig:sfig3}
\end{subfigure}
\caption{ Different elimination process on different action} \label{fig:E-3}
\label{fig:fig}
\end{figure}

Let's see figure \ref{fig:E-3}. A tiny matrix with only 3 rows. If only row transformation is considered, there are six different elimination process in total. 
\begin{equation}
    \left\{\begin{matrix} 
    R_{1}\rightarrow R_{2}\rightarrow R_{3}
    \\ R_{1}\rightarrow R_{3}\rightarrow R_{2}
    \\ R_{2}\rightarrow R_{1}\rightarrow R_{3}
    \\ R_{2}\rightarrow R_{3}\rightarrow R_{1}
    \\ R_{3}\rightarrow R_{1}\rightarrow R_{2}
    \\ R_{3}\rightarrow R_{2}\rightarrow R_{1}
\end{matrix}\right.
\end{equation}

Figure \ref{fig:E-3} shows three case, which is on the different reward and action.
\begin{itemize}
\item  Figure 2.(a) $R_{1}\rightarrow R_{2}\rightarrow R_{3}$. No pivoting. No row permutation. This process requires minimal computation. 
\item  Figure 2.(b) $\left [ R_1,R_2 \right ] \rightarrow R_{3}$.  $R_{1}, R_{2}$ can be eliminated at the same time, so there is a high degree of parallelism.
\item  Figure 2.(c) $R_{3}\rightarrow R_{2}\rightarrow R_{1}$. Apply column pivoting strategy to ensure numerical stability.
\end{itemize} 
So even for this tiny matrix, there are different strategy for different goals. At each step, there are different policy to pick different action. We could use off-line Q-Learning to learn these strategies. Then apply it these strategies to improve the efficiency and accuracy of the solution process.

\subsection{GSS - a High Performance Solver with Q-Learning based Scheduler}

We have implemented Q-Learning based scheduling algorithm in GSS\cite{Chen2018}. It use Q-Learning based task-load-tree split technique to create a sub-matrix that can be fully factored in GPU. So the data transfer time reduced to a minimum. GSS also uses some adaptive load balance techniques in the framework of Q-Learning. In the static scheduling stage, LU factorization is split into many parallel tasks. Each task is assigned to a different computing core (CPU, GPU, ...). Then in the actual elimination process, the dynamic scheduling strategy will do more load balancing on pretrained $Q$ table. That is, if GPU cores have high computing power, then it will run more tasks automatically. If CPU is more powerful, then GSS will give it more tasks. The choice of target core is learned from $Q$ table. 

There are some experimental results in https://github.com/closest-git/GSS. For many large matrices, GSS is about 2-3 times faster than PARDISO\cite{schenk2004solving} in MKL. We would give more detailed results in subsequent papers.

\section{Prospect}
In this paper, we list many combinatorial optimization processes in sparse solver. We unified these processes into the framework of Markov Decision Process, with detailed analysis of state space, action space, and reward. Then use Q-Learning technique to improve the solution. 

This is only the first step to rethinking and recasting the Gaussian elimination process. More work should be done. We will report more method and numerical experimental results in subsequent papers.

\bibliographystyle{unsrt}  
\bibliography{references}  






\end{document}